\documentclass[12pt]{article}
\usepackage{amsmath}
\usepackage{amsfonts}
\usepackage{amssymb}
\usepackage{graphicx}
\usepackage{color}
\usepackage{amsthm}
\usepackage{hyperref}

\theoremstyle{plain}

\def\DHLhksqrt#1#2{\setbox0=\hbox{$#1\sqrt{#2\,}$}\dimen0=\ht0
\advance\dimen0-0.2\ht0
\setbox2=\hbox{\vrule height\ht0 depth -\dimen0}%
{\box0\lower0.4pt\box2}}
\def\Re{\mathfrak{Re}}
\def\Im{\mathfrak{Im}}

\def\bee{\begin{equation}}
\def\eee{\end{equation}}

\def\DHLhksqrt#1#2{\setbox0=\hbox{$#1\sqrt{#2\,}$}\dimen0=\ht0
\advance\dimen0-0.2\ht0
\setbox2=\hbox{\vrule height\ht0 depth -\dimen0}%
{\box0\lower0.4pt\box2}}

\allowdisplaybreaks

\begin{document}

\bigskip
\bigskip
\bigskip
\bigskip
\bigskip
\bigskip
\centerline{\Large\bf Almost  all of the nontrivial  zeros   }
\bigskip
\centerline{\Large\bf   of the Riemann zeta-function  are on the critical line}

\bigskip
\bigskip

\begin{center}
{\large \sl{  \sl{ C. Dumitrescu}$^1$, \sl{M. Wolf }$^2$} }  
\end{center}

\begin{center}
$^1$Kitchener,  Canada, email:  cristiand43@gmail.com\\
$^2$Cardinal  Stefan  Wyszynski  University, Warsaw,   Poland, e-mail:  primes7@o2.pl
\end{center}

\bigskip
\bigskip

\begin{center}
{\bf Abstract}\\
\bigskip
\begin{minipage}{12.8cm}
Applying   Littlewood's lemma in connection to Riemann's Hypothesis  and  exploiting the symmetry of Riemann's $xi$ function
we show that almost all  nontrivial Riemann's Zeta zeros are on the critical line.
\end{minipage}
\end{center}

\bigskip\bigskip

\section{Introduction}

In his only paper devoted to the number theory published in 1859 \cite{Riemann-1859} 
(it was also included as an appendix in \cite{Edwards})  Bernhard Riemann  continued  analytically  the series
\bee
\sum_{n=1}^\infty \frac{1}{n^s}, ~~~~~~s=\sigma+it,   ~~~\sigma>1
\eee
to the complex plane  with exception of $s=1$, where the above series is a harmonic divergent series. He has done it using the integral
\bee
\zeta(s)=\frac{\Gamma(1-s)}{2\pi i}\int_\mathcal{C} \frac{(-z)^s}{e^z-1} \frac{dz}{z}, 
\label{zeta-plane}
\eee
\noindent
\noindent where  the contour $\mathcal{C}$ is \\

\centerline{      }

\centerline{      }

\centerline{      }

\centerline{      }

\centerline{      }

\centerline{      }

\begin{picture}(0,0)(0,0)
\thicklines
\put(200,30){\vector(0,1){70}}
\put(140,60){\vector(1,0){140}}
\put(230,70){\mbox{$\mathcal{C}$}}
\put(200,60){\oval(40,40)[l]}
\put(200,65){\oval(30,30)[rt]}
\put(200,55){\oval(30,30)[rb]}

\put(215,55){\vector( 1,0){45}}
\put(260,65){\vector(-1,0){45}}
\end{picture}

\vskip-1cm

\noindent  The definition of $(-z)^s$ is $(-z)^s= e^{s\log (-z)}$, where
the definition of $\log(-z)$ conforms to the usual definition of $\log(z)$ for $z$ not
on the negative real axis as the branch which is real for positive real $z$, thus
$(-z)^s$ is not defined on the positive real axis, see \cite[p.10]{Edwards}.  Appearing in  \eqref{zeta-plane}  the  gamma function  $\Gamma(z)$
has many representations, we present  the Weierstrass product:
\bee
\Gamma(z) = \frac{e^{-C z}}{z} \prod_{k=1}^\infty \frac{ e^{z/k}}{\left(1 + \frac{z}{k}\right)} ~.
\label{Gamma}
\eee
Here  $C$ is the Euler--Mascheroni constant
\bee
C = \lim_{n \rightarrow \infty}  \left( \sum_{k=1}^n \frac{1}{k} - \log(n) \right) = 0.577216\ldots.
\eee
From (\ref{Gamma}) it is seen that $\Gamma(z)$  is defined for all complex numbers  $z$, except $z = -n$ for integer  $n > 0$,
where are the  simple poles of $\Gamma(z)$.  The most popular definition of the gamma function given by the integral
$\Gamma(z) = \int_0^\infty e^{-t} t^{z-1} dt $  is valid only for $\Re[z]>0$.
Recently perhaps  over 100  representations of  $\zeta(s)$ are known,
for review of the integral and  series  representations see \cite{Milgram-2013}.

The function  $\zeta(s)$ has two kinds of zeros:  trivial zeros  at $s=-2n, ~n=1, 2, 3, \ldots$
and nontrivial zeros  in  the critical strip $0<\Re(s)<1$.  In  \cite{Riemann-1859}   Riemann  made the assumption,  now called the
{\it  Riemann Hypothesis} (RH  for short in following), that all nontrivial zeros $\rho_n$ lie on the {\it critical line} $\Re[s]=\frac{1}{2}$:
$\rho_n=\frac12+i\gamma_n$. Contemporary the above requirement is augmented by the demand that all nontrivial zeros are simple.
Riemann has shown that $\zeta(s)$ fulfills the {\it  functional identity}:
\bee
\pi^{-\frac{s}{2}}\Gamma\left(\frac{s}{2}\right)\zeta(s)
           =\pi^{-\frac{1-s}{2}}\Gamma\left(\frac{1-s}{2}\right)\zeta(1-s),~~~~{\rm  for } ~ s \in\mathbb{C}\setminus\{0,1\}.
\label{functional-zeta}
\eee
The above form of the functional equation is explicitly symmetrical with respect to the line $\Re(s)=1/2$: the change
$s\rightarrow 1 - s$ on both sides of (\ref{functional-zeta}) shows that the values of the combination
of functions $\pi^{-\frac{s}{2}}\Gamma\left(\frac{s}{2}\right)\zeta(s)$  are  the same at points  $s$  and
$s-1$.  Thus it is convenient to introduce the function
\bee
\xi(s)=\frac{1}{2}s(s-1)\Gamma\left(\frac{s}{2}\right)\zeta(s).
\label{def-xi}
\eee
Then  the functional identity takes the simple form:
\bee
\xi(1-s)=\xi(s)
\label{funct-eq}
\eee

\begin{figure*}[ht]
\includegraphics[page=1,width=.95\textwidth]{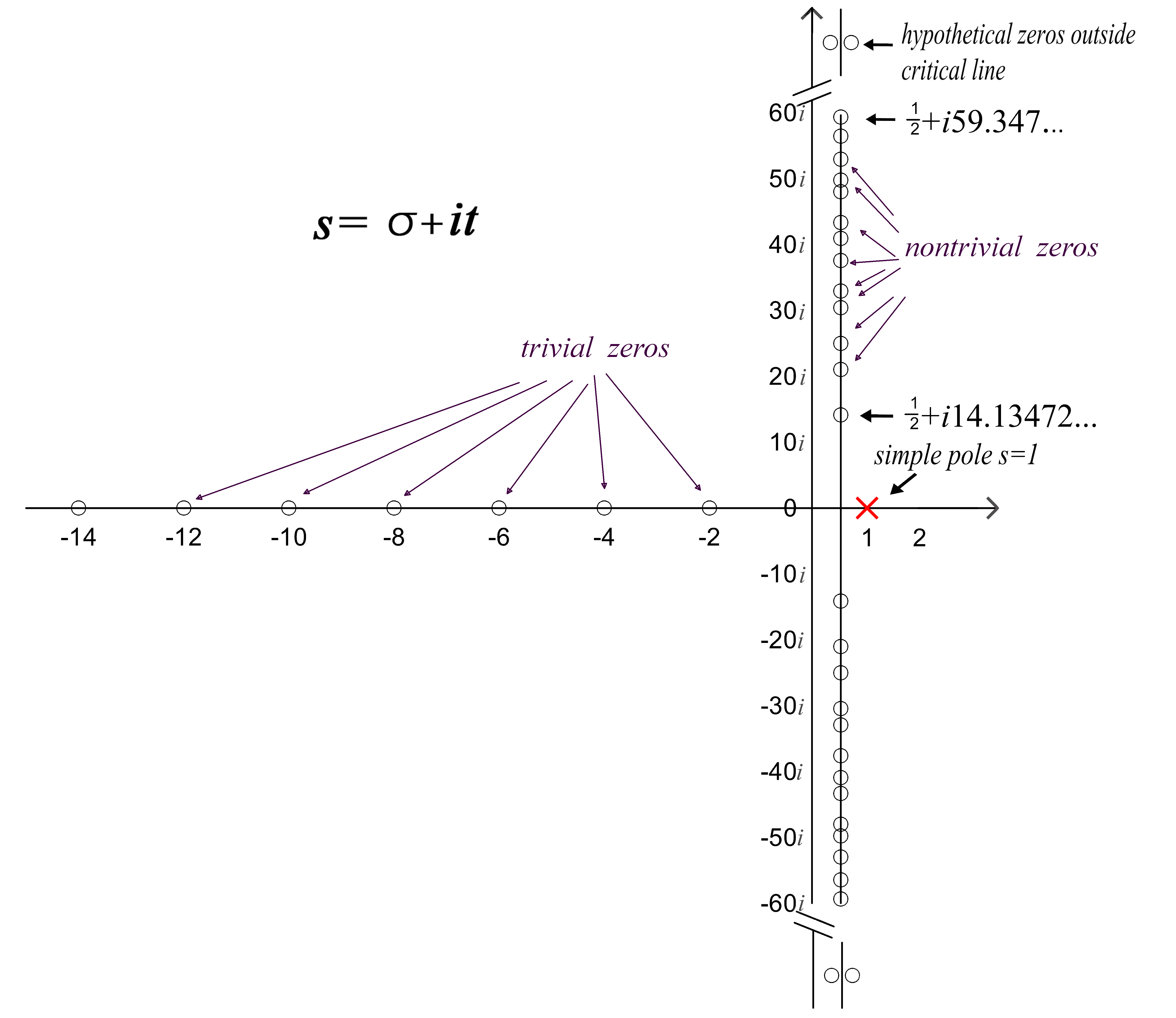}
\caption{\small The location of zeros of the Riemann $\zeta(s)$ function.}
\bigskip
\label{plane}
\end{figure*}

The fact that $\zeta(s)\neq 0$ for $\Re(s)>1$ and the  form  of the  functional
identity entails that {\it nontrivial}  zeros $\rho_n=\beta_n + i\gamma_n$ are located in the {\it critical  strip}:
\[
0 \leq  \Re [\rho_n] =\beta_n \leq 1.
\]
From the complex  conjugation of $\zeta(s)=0$  it follows that if  $\rho_n=\beta_n + i\gamma_n$ is a zero, then
$\overline{\rho_n}=\beta_n -i\gamma_n $ also is a zero.   From the symmetry of the functional equation (\ref{functional-zeta}) with
respect to  the line $\Re[s]=\frac{1}{2}$  it follows, that if $\rho_n=\beta_n + i\gamma_n$ is a zero, then  $1-\rho_n=1-\beta_n - i\gamma_n$  and     $
1-\overline{\rho_n}=1-\beta_n + i\gamma_n$  are also zeros: they are located symmetrically  around  the straight line  $\Re[s]=\frac{1}{2}$ and
the axis $t=0$,  see Fig. \ref{plane}.

The classical (from XX  century)  references on the RH are \cite{Titchmarsh}, \cite{Edwards},  \cite{Ivic1985}, \cite{Karatsuba-Voronin}.
In the XXI century there appeared two monographs about the zeta function: \cite{Borwein_RH}  and    \cite{Broughan_2017}.

There was a  lot of attempts to prove RH and the common  opinion  was that it is  true.  However let us notice that there were  famous
mathematicians: J. E. Littlewood  \cite[p.345]{Burkill}, \cite[p. 390]{Good1962},    P. Turan and A.M. Turing  \cite[p.1209]{Booker_turingand}
M. Huxley \cite[p. 357]{Derbyshire}  believing  that the  RH is not true,   see  also  the paper ``On some reasons for doubting the
Riemann  hypothesis'' \cite{Ivic2003}  (reprinted in \cite[p.137]{Borwein_RH})  written by  A. Ivi{\'c}.
In Karatsuba's talk  \cite{linki_Karatsuba} at 1:01:10 \footnote{We thank A. Kourbatov for bringing  this fact to our  attention}
he mentions that   Atle Selberg had serious doubts whether RH is true or not.   New arguments against RH   can be found in \cite{Blanc_2019}.
When  J. Derbyshire asked A. Odlyzko about his opinion on the validity of RH   he replied ``Either it's true, or else it isn't''
\cite[p. 357--358]{Derbyshire}. There were  some attempts to prove RH using the  physical  methods,  see  \cite{Physics-and-RH-RevModPhys}
or  \cite{Wolf_2020_RPP}.

In \cite[p.81]{Steinhaus_2015} we read: ``Hilbert said that if he could rise
from the dead in 200 years, his first thought would not be to ask what social or
technological progress there had been, but what had been discovered about the zeros
of the zeta function  $\zeta$ because that is not only the most interesting unanswered
mathematical question, but the most interesting of all questions...''

\section{Zeta's  zeros on the critical  line}

In  1914  G.H. Hardy \cite{Hardy1914}  (reprinted in \cite{Borwein_RH})  proved first  result in favor  of the RH:  there are infinitely many
zeros of $\zeta(s)$  on the critical line. Let $N_0(T)$  denotes  number of the $\zeta(s)$  zeros  on the critical line $\frac12+i\gamma_n$ with imaginary
part $0<\gamma_n<T$.  In 1921   Hardy and  Littlewood in   \cite{Hardy_L-1921}   proved that $N_0(T)>const\cdot T$  for  large  $T$.
Let $N(T)$ denotes the  number   of zeta zeros  $\rho_n=\beta+i\gamma_n$  in the critical  strip up to  $T$, i.e. the number
of zeta zeros in the rectangle $0<\beta<1, ~~0<\gamma_n<T$.    In 1942 A. Selberg  in \cite{Selberg1942b}
proved that  $N_0(T)>const \cdot T\log  \log \log(T)  $ for  large  $T$  and in \cite{Selberg1942}  he improved it to $N_0(T)>const \cdot T \log(T)  $
for  large  $T$ (these  two papers are reprinted  in \cite{Selberg01}).    Norman Levinson  proved  \cite{Levinson1974} that more than one-third of
zeros  of  Riemann's $\zeta(s)$ are  on critical line $N_0(T)>N(T)/3$ by relating
the zeros of the zeta function to those of its derivative.  Later Levinson  \cite{Levinson_1975}  improved the proportion of zeros on the  critical
line to   $0.3474$.  Brian   Conrey    in  1989  improved this further to two-fifths  (precisely 40.77 \%)  \cite{Conrey1989}.  Next
S. Feng  proved that at least 41.28\% of the zeros of the Riemann zeta  function are on the   critical line \cite{Feng}.
The present record  seems to belong to K.  Pratt et.  al. who proved that at least $5/12=0.41666\ldots$ of the zeros of the Riemann zeta
function are on the   critical line \cite{Pratt_2019}.

In this paper we are going to apply the  Littlewood's Lemma to  show that almost all zeta zeros are on the critical line:
\bigskip

{\bf Theorem:}  Almost  all zeros $\rho_n=\beta_n + i \gamma_n$ of the $\zeta(s)$  function  have $\beta_n=\frac12$.  Here by ``almost all'' we mean, that
\bee
\frac{N_0(T)}{N(T)}=1 + \mathcal{O}\left(\frac{1}{T\log(T)}\right),
\label{tw}
\eee
where  $T$ is not an imaginary part of the nontrivial  zeta's  zero.

\bigskip
\bigskip

\section{Proof of the Theorem}.

We will use the Littlewood's Lemma  (see e.g.  \cite[Chap.21]{Iwaniec-2014}:

\bigskip

{\bf Littlewood's Lemma:}   ~Let  $F(s)$ be  holomorphic   function inside rectangle $\mathcal{D} $   with sides  parallel  to  axes not vanishing  on the
boundary $\partial \mathcal{D}$  and  let  dist$(\rho)$  denotes  the  distance of the zero $\rho $ of $F(s)$,  i.e.  $F(\rho)=0$,
to the left  side of  $\mathcal{D}$.   Then  
\bee
\sum_{\rho \in \mathcal{D}} {\rm dist}(\rho)=\frac{1}{2\pi}\oint_{\partial \mathcal{D}}\log(F(s)) ds,  ~~~~s=\sigma+it
\label{LL}
\eee

For the function $F(s)$  we will substitute  the  Riemann's $\xi(s)$  function  defined in \eqref{def-xi}  which has  the same zeros as
$ \zeta(s)$:
\bee
\xi(s)=\frac12 s(s-1)\pi^{-s/2}\Gamma\left(\frac{s}{2}\right) \zeta(s).
\label{ksi}
\eee
In  addition to functional equation \eqref{funct-eq} it fulfils:

\bee
\overline{\xi(s)}=\xi(\overline{s}).
\label{symm2}
\eee
We have
\bee
\log(\xi(s))=\log|\xi(s)|+i\arg (\xi(s)).
\eee
Let $0<\alpha<\frac12$ be a  real  number. We consider the rectangle $\mathcal{D}(\alpha, T) $  with  vertices  $A=(1-\alpha, iT), B=(\alpha, iT),
C=(\alpha, -iT) , D=(1-\alpha, -iT)$, see Fig.2  and  we will look for zeros of $\xi(s)$  which
are  the same as zeros of $\zeta(s)$.  We have in  this case  (see .eq(21) in \cite{Iwaniec-2014}):
\[
\sum_{\rho \in \mathcal{D(\alpha)}} {\rm dist} \rho=\frac{1}{2\pi}\int_{-T}^T\big(\log(|\xi(\alpha+it)|)-\log(|\xi(1-\alpha+it)|)\big) dt
\]
\bee
+\frac{1}{2\pi}\int_{\alpha}^{1-\alpha}\big(\arg(\xi(\sigma+iT))-\arg(\xi(\sigma-iT)) \big)d\sigma,
\label{Calki}
\eee
where the argument is defined by continuous variation starting with any fixed value  at a chosen point.
  .

From   functional identity  \eqref{funct-eq}  and \eqref{symm2}  we have
\bee
|\xi(\alpha+it)|=|\xi(1-\alpha-it)|=|\xi(1-\alpha+it)|
\eee
thus the first integral in \eqref{Calki} vanishes. Next, we know that
\bee
\arg(\xi(\sigma-iT))=-\arg(\xi(\sigma+iT))
\eee
and it follows
\bee
\sum_{\rho \in \mathcal{D(\alpha)}} {\rm dist}(\rho)=\frac{1}{\pi}\int_\alpha^{1-\alpha}\arg(\xi(\sigma+iT)) d\sigma
\label{rhs}
\eee

First we will calculate rhs of \eqref{rhs}.    We have:
\bee
\arg(\xi(s))=\arg\left(\frac12 s(s-1)\right)+\arg\left(\pi^{-s/2}\right)+\arg\left(\Gamma\left(\frac{s}{2}\right)\right)+\arg(\zeta(s))
\eee

From the proviso  in \eqref{Stirling}  we have to restrict $s$ to  principal   branch,  thus  argument  for  $s=\sigma+iT$   is  close to
$\pi/2$   because $\arg(s)$ is a little bit less than  $\pi/2$  and $\arg(s-1)$  is a little bit more  than  $\pi/2$  thus  we write
\bee
\arg(s(s-1))=\pi +o(1).
\eee
Hence we have
\bee
\int_\alpha^{1-\alpha} \arg\left(\frac12 s(s-1)\right) d\sigma = (1-2\alpha)\left(\pi+ o(1)\right).
\label{calka_stalej}
\eee
From Stirling formula,  see  e.g.  \cite[eq. 6.1.37 and eq.(6.1.41)]{abramowitz+stegun}   we have
\bee
\log(\Gamma(z))=\left(z-\frac12\right)\log(z)-z+\frac12\log(2\pi)+\frac{1}{12z}+\ldots, ~~~~~\arg(z)<\pi
\label{Stirling}
\eee
and in our case:
\bee
\log\left(\Gamma\left(\frac{s}{2}\right)\right)=\left(\frac{\sigma}{2}-\frac12+i\frac{T}{2}\right)\left(\log\left(\frac{\sigma^2}{4}+\frac{T^2}{4}\right)^\frac12+i\arg\left(\frac{\sigma}{2}+i\frac{T}{2}\right)\right)
\eee
\[
-\left(\frac{\sigma}{2}+i\frac{T}{2}\right)+\frac12 \log(2\pi)+o(1)
\]
Hence
\bee
\arg\left(\Gamma\left(\frac{s}{2}\right)\right)=\Im(\log\left(\Gamma\left(\frac{s}{2}\right)\right)=
\eee
\[
\frac{T}{4}\log\left(\frac{\sigma^2}{4}+\frac{T^2}{4}\right)-\frac{T}{2}+\left(\frac{\sigma}{2}-\frac12\right)\arg\left(\frac{\sigma}{2}+i\frac{T}{2}\right)
\]
For the last term  we can write for  very large $T$
\bee
\left(\frac{\sigma}{2}-\frac12 \right)\arg\left(\frac{\sigma}{2}+i\frac{T}{2}\right)=\left(\frac{\sigma}{2}-\frac12 \right)\left(\frac{\pi}{2}+o(1)\right)
\eee
For very large $T$  we can skip terms with $\alpha$ under logarithms:
\bee
T\log\left(\frac{\alpha^2}{4}+\frac{T^2}{4}\right)=T\log\left(\frac{T^2}{4}\left(1+\frac{\alpha^2}{T^2}\right)\right)=
\eee
\[
2T\log\left(\frac{T}{2}\right)+T\log\left(1+\frac{\alpha^2}{T^2}\right)=2T\log\left(\frac{T}{2}\right) + \mathcal{O}\left(\frac{1}{T}\right).
\]
Further  we have
\bee
\arg(\pi^{-(\sigma+iT)/2})=-\frac{T}{2}\log(\pi)
\label{Pi}
\eee
Together  from  above equations  we have:
\bee
\arg(\xi(\sigma+iT))=\frac{T}{2}\log\left(\frac{\sigma^2}{4}+\frac{T^2}{4}\right)+C(\sigma)-\frac{T}{2}(1+\log(\pi))+\arg(\zeta(\sigma+iT))+o(1).
\eee
Above $C(\sigma)$  absorbs  constants:
\bee
C(\sigma)=\pi+\frac{\pi}{2}\left(\frac{\sigma}{2}-\frac12 \right) +o(1)=\frac{3\pi}{4}+\frac{\pi}{4}\sigma+o(1)
\eee

\setlength{\unitlength}{1.5cm}
\begin{picture}(10,10)(-5,-4)
\thicklines
\put(-2.5,0){\vector(1,0){5}}
\put(2.7,-0.1){$\sigma$}
\put(2.2,0.1){$1$}
\put(2.1,-0.1){\line(0,1){0.2}}
\put(-1.5,-4.5){\vector(0,2){9}}
\put(0.3,-4.2){\color{red}\line(0,1){8.4}}
\put(0.5,4){${\small \frac12 +it}$}
\put(0.3,3.0){${\small \frac12 +iT}$}
\put(0.6, 3.2){\vector(-1,1){0.3}}

\multiput(1.3,-4)(0,2.5){3}{\vector(0,1){2.5}}
\multiput(-0.7,-1.5)(0,2.5){3}{\vector(0,-1){2.5}}
\put(1.3,3.5){\vector(-1,0){2}}
\put(-0.7,-4){\vector(1,0){2}}
\put(-1.3,1.5){\vector(-1,0){0.2}}
\put(-1.2,1.4){${\small\displaystyle{\small{\alpha}}}$}
\put(-0.9,1.5){\vector(1,0){0.2}}

\multiput(2.1,-4)(0,0.8){11} {\line(0,1){0.2}}

\put(1.3,3.5){$A$}
\put(-1,3.5){$B$}
\put(-1,-4.2){$C$}
\put(1.3,-4.2){$D$}

\put(-1.8,4.5){$t$}
\put(-4, -4.8){Figure 2. ~Rectangle $\mathcal{D}(\alpha, T)$,  in red critical line is plotted.}
\put(-4.8, -5.2) { Sides  $AB$ and $CD$  have length  $1-2\alpha$  while  sides $CB$ and $DA$  have length $2T$. }\\
\end{picture}

\bigskip
\bigskip

\centerline{      }

\centerline{      }

\centerline{      }

\bigskip
\bigskip

Integrating by parts   we obtain
\bee
\int\log\left(\frac{\sigma^2}{4}+\frac{T^2}{4}\right)d\sigma=\sigma\log\left(\frac{\sigma^2}{4}+\frac{T^2}{4}\right)-2\sigma+2T\arctan\left(\frac{\sigma}{T}\right)  +constant
\eee
and it gives
\bee
\int_\alpha^{1-\alpha}\arg(\xi(\sigma+iT)) d\sigma=(1-\alpha)\frac{T}{4}\left(\log\left(\frac{(1-\alpha)^2}{4}+\frac{T^2}{4}\right)-2\right)
\eee
\[
-\alpha\frac{T}{4}\Big(\log\left(\frac{\alpha^2}{4}+\frac{T^2}{4}\right)-2\Big)+\frac{T^2}{2}\Big(\arctan\left(\frac{1-\alpha}{T}\right)-\arctan\left(\frac{\alpha}{T}\right)\Big)
\]
\[
-(1-2\alpha)\frac{T}{2}(1+\log(\pi))+C_2+\int_\alpha^{1-\alpha}\arg(\zeta(\sigma+iT))d\sigma+\mathcal{O}(1)+C_2
\]
where
\bee
C_2=\int_\alpha^{1-\alpha} C(\sigma) d\sigma=(1-2\alpha)\left(\frac{7 \pi}{8}+o(1)\right)
\eee

Using the Maclaurin series  expansion  of $\arctan(x)$  we find  that for large $T$
\bee
T\Big(\arctan\left(\frac{1-\alpha}{T}\right)-\arctan\left(\frac{\alpha}{T}\right)\Big)=1-2\alpha+\mathcal{O}\left(\frac{1}{T}\right).
\eee
Adding  \eqref{calka_stalej}  and  \eqref{Pi} gives us finally  that the  rhs of \eqref{rhs}  is equal to
\[
\frac{1}{\pi}\int_\alpha^{1-\alpha} \Big\{\arg\left(\frac12 s(s-1)\right)+\arg\left(\pi^{-s/2}\right)+\arg\left(\Gamma\left(\frac{s}{2}\right)\right)+\arg(\zeta(s)\Big\}  d\sigma
\]
\bee
=(1-2\alpha)\Big(\frac{T}{2\pi}\log\left(\frac{T}{2\pi}\right)-\frac{T}{2\pi}+\frac{C_2}{\pi}+o(1)\Big)+\frac{1}{\pi}\int_\alpha^{1-\alpha} \arg(\zeta(s)\Big\}  d\sigma .
\label{rhs2}
\eee

We give estimate for the integral of the argument of $\zeta(s)$.  The  mean value theorem for definite integrals:  Let $f : [a, b] \to  \mathbb{R}$
be a continuous function. Then there exists $c \in (a, b)$ such that

\bee
 \int _{a}^{b}f(x)\,dx=f(c)(b-a).
\eee
Applying it to $\arg(\zeta(\sigma+iT))$, when $T$ is  not equal to imaginary  part of nontrivial zero ($\arg(\zeta(s))$ is not
defined  for $s$  equal to  zero of zeta) we have:
\bee
\int_\alpha^{1-\alpha} \arg(\zeta(\sigma+iT)) d\sigma = (1-2\alpha) \arg(\zeta(\sigma_0+iT))
\eee
for some $\sigma_0 \in (\alpha, 1-\alpha)$.  Letting $\alpha\to \frac12 $ we have
\bee
\lim_{\alpha\to \frac12} \frac{1}{1-2\alpha}\int_\alpha^{1-\alpha} \arg(\zeta(\sigma+iT)) d\sigma = \arg\left(\zeta\left(\frac12 +iT\right)\right)
\eee

Thus finally  we  obtain  for  $\alpha$  close  to $\frac12$
\[
\frac{1}{\pi}\int_\alpha^{1-\alpha}\arg(\xi(\sigma+iT))d\sigma  =
\]
\bee
=(1-2\alpha)\Big\{\frac{T}{2\pi}\log\left(\frac{T}{2\pi}\right)-\frac{T}{2\pi}+\frac78 +\frac{1}{\pi}\arg\left(\zeta\left(\frac12 +iT\right)\right)+o(1)\Big\}.
\label{rhs3}
\eee
From the Riemann--von Mangoldt formula we have  that the number of zeta zeros $N(T)$  with positive imaginary parts $<T$ (and real part inside critical
strip)  is  (see eq.(2.3.6)  in  \cite{Borwein_RH})
 \bee
 N(T)=\frac{T}{2\pi}\log\left(\frac{T}{2\pi}\right)-\frac{T}{2\pi}+\frac78+\frac{1}{\pi}\arg\zeta\left(\frac12+iT\right) + \mathcal{O}(T^{-1}).
 \label{Mangoldt}
 \eee
 The term  $7/8$ was already  known  to von Mangoldt, see \cite[ p.10  eq.(8)]{Mangoldt1900}.  Comparing \eqref{rhs3}  and  \eqref{Mangoldt}
 we obtain
\bee
\frac{1}{\pi}\int_\alpha^{1-\alpha}\arg(\xi(\sigma+iT))d\sigma  = (1-2\alpha)\left(N(T)+o(1)\right)
\label{main}
\eee

\bigskip

Now we will  focus on the lhs of  \eqref{rhs}. This needs  information about specific  location of  zeros. Because there are only finite number
of zeta zeros inside rectangle $ABCD$,  we can choose $\alpha$ so close
to  $\frac12 $ that zeros  off critical line be outside rectangle  $\mathcal{D}(\alpha, T) $, i.e.  our symmetric rectangle
closely hugs the critical   line that inside  are   only zeros  on $\Re(s)=\frac12$.  Then we have
\bee
 \sum_{\rho \in \mathcal{D}(\alpha, T)} {\rm dist}(\rho)=\left(\frac12-\alpha\right)2N_0(T).
 \label{lhs2}
\eee
Comparing \eqref{main} and \eqref{lhs2}  we can cancel  $1-2\alpha\neq 0$ on both sides  what gives:
\bee
N_0(T)=N(T)+o(1)
\eee
and we  obtain  \eqref{tw}.  \hfill$\Box~~~~~~~~$

\bigskip

\bigskip

\section{Comments}

The von  Mangoldt formula   \eqref{Mangoldt}   was  obtained from the Argument   Principle.   This principle  gives the total number of zeros
in some region   and it  needs   the information  on the change of the  argument.  The Littlewood  lemma is more  powerful:  it involves
specific  locations  of zeros and the integrals  of  the argument.

Our method is based on the symmetry relations satisfied by the $\xi(s)$  function (in relation to Littlewood's lemma),  which drastically simplify
the calculations.  The other two major ideas in our work are, first, the use of the mean value theorems for integrals, (and second)  in
combination with that limit process when the rectangle hugs the critical line together with the behavior of various terms in our calculations
under this process. So there are about three fundamental ideas that we use in our work.  We rely more on symmetry and the consequences of
a limit process,  rather than the computational techniques involving smoothing functions (mollifiers),  like Conrey and Levinson employ.
When we apply Littlewood's lemma,  we know that there exists an analytic branch of $\log(\xi)$ on the simply connected domain represented
by our rectangle without the horizontal branch cuts connecting the left side of our rectangle with  $\zeta$  zeros inside our rectangle.
For any branch of the logarithm,  in other words for any branch of the argument, the conclusions of our calculations are not
significantly affected.

\bigskip

Because   in \eqref{tw}  in the limit $T\to \infty $  there will be equality $N_0(T)=N(T)$  we  propose

\bigskip

{\bf Conjecture:}  There exists  such a constant $Y>0$ that if $\beta +i\gamma$  is  a zero of $\zeta(s)$  and $|\gamma|>Y$   then $\beta=\frac12$.

\bigskip

If  RH is true then the above Conjecture is empty,  i.e.  $Y=0$.

In the paper  \cite{Churchhouse1968}    I.J. Good and R.F. Churchhouse  gave the arguments that RH  is true with probability  1.
Above we have showed that  almost 100\%   of the $\zeta(s)$  are on the critical line.

\end{document}